\documentclass{amsart}

\usepackage{amsmath,amssymb,amsthm}

\newtheorem{prop}{Proposition}
\newtheorem{theorem}[prop]{Theorem}
\newtheorem{lemma}[prop]{Lemma}

\theoremstyle{definition}
\newtheorem{definition}[prop]{Definition}

\DeclareMathOperator{\lt}{lt}
\DeclareMathOperator{\supp}{supp}

\newcommand{\inv}{^{-1}}
\newcommand{\sidehat}{^{\wedge}}
\newcommand{\doubledual}[1]{\,\widehat{\!\widehat{#1}}}

\allowdisplaybreaks
\title{Group Bundle Duality}
\author{Geoff Goehle}
\address{Department of Mathematics, Dartmouth College 
        Hanover, NH 03755 }
\email{goehle@dartmouth.edu}

\begin{document}

\begin{abstract}
This paper introduces a generalization of 
Pontryagin duality for locally compact Hausdorff Abelian groups
to locally compact Hausdorff Abelian group bundles. 
\end{abstract}

\maketitle

First, recall that a group bundle is just a groupoid where the 
range and source maps coincide.  An Abelian group bundle is 
a bundle where each fibre is an Abelian group.  When working with a
group bundle $G,$ we will use $X$ to denote the unit space of $G$ and 
$p:G\rightarrow X$ to denote the combined range and source maps.  Furthermore,
we will use $G_x$ to denote the fibre over $x$. 
Group bundles, like general groupoids, may not have a Haar system 
but when they do the Haar system has a special form.  If $G$ is a locally
compact Hausdorff group bundle with Haar system, denoted by
$\{\beta^x\}$ throughout the paper, then 
$\beta^x$ is Haar measure on the fibre $G_x$ for all $x\in X$.  At
this point, it is convenient to make the standing assumption that all
of the locally compact spaces in this paper are Hausdorff.  

Now suppose $G$ is an Abelian, second countable, locally compact
group bundle with Haar system $\{\beta^x\}$. 
Then $C^*(G,\beta)$ is a separable Abelian $C^*$-algebra and 
in particular $\widehat{G}=C^*(G,\beta)\sidehat$ is a second countable 
locally compact Hausdorff space \cite[Theorem 1.1.1]{invitation}.  We cite 
\cite[Section 3]{ctgIII}
to see that each element of $\widehat{G}$ is of the form $(\omega,x)$ with 
$x\in X$ and $\omega$ a character in the Pontryagin dual of $G_x$, 
denoted $(G_x)\sidehat$.  The
action of $(\omega,x)$ on $C_c(G)$ is given by 
\begin{equation}
\label{eq:1}
(\omega,x)(f) = \int_G f(s)\omega(s) d\beta^x(s).
\end{equation}
Since every element in $\widehat{G}$ is a character on a fibre of 
$G,$ we are justified in thinking of $\widehat{G}$ as a bundle over 
$X$ with fibres $\widehat{G}_x = (G_x)\sidehat$ and action on $C^*(G,\beta)$ 
given by \eqref{eq:1}.  We will use $\hat{p}$ to denote the projection from 
$\widehat{G}$ to $X$ and $\omega$ to denote the 
element $(\omega,\hat{p}(\omega))$ in $\widehat{G}$.  

At this point, it is clear that $\widehat{G}$ is algebraically a 
group bundle.  In order for it to be a topological groupoid, 
 we must show that the groupoid operations are continuous with 
respect to the Gelfand topology on $\widehat{G}$.  To this end, we reference
the following characterization of the topology on $\widehat{G}$.  
\begin{lemma}[{\cite[Proposition 3.3]{ctgIII}}]
\label{prop:1}
Let $G$ be a second countable locally compact abelian group bundle 
with Haar system. Then a sequence $\{\omega_n\}$ in 
$\widehat{G}$ converges to $\omega_0$ in $\widehat{G}$ if and only if 
\begin{enumerate}
\item $\hat{p}(\omega_n)$ converges to $\hat{p}(\omega_0)$ in $X$, and 
\item if $s_n\in G_{\hat{p}(\omega_n)}$ for all $n\geq 0$ and $s_n\rightarrow s_0$ in $G$, then 
$\omega_n(s_n)\rightarrow \omega_0(s_0)$.  
\end{enumerate}
\end{lemma}
The first thing we can conclude from this lemma is that the 
restriction of the topology on $\widehat{G}$ to $\widehat{G}_x$ is the same
as the topology on $\widehat{G}_x$ as the dual group of $G_x$.  The 
second thing we conclude is that the topology on $\widehat{G}$ is independent 
of the Haar system $\beta$.  Furthermore, recall that
the groupoid operations on $\widehat{G}$ are those coming from the 
dual operations on $\widehat{G}_x$.  In other words, the operations are pointwise
multiplication and conjugation of characters, and it follows from 
Lemma \ref{prop:1} that these operations are continuous.  Therefore,
we have proven the lemma. 
\begin{lemma}[{\cite[Corollary 3.4]{ctgIII}}]
\label{thm:dual}
Let $G$ be a second countable locally compact Abelian group bundle with 
Haar system.  Then $\widehat{G}$, equipped with the Gelfand 
topology, is a second countable locally compact Abelian group 
bundle with fibres $\widehat{G}_x = (G_x)\sidehat$.  
\end{lemma}
Now we can make our first definition. 
\begin{definition}
If $G$ is a second countable locally compact Abelian group bundle 
with Haar system, then we define the dual bundle to be 
$\widehat{G}=C^*(G)\sidehat$ equipped with the groupoid operations coming 
from the identification of $\widehat{G}_x$ as the dual of $G_x$.  
We will use $\hat{p}$ to denote the projection on 
this bundle.
\end{definition}

This definition gives rise to the notion of a duality theorem for
group bundles.  The main result of this paper is to 
prove the following theorem, stated without proof in
\cite[Proposition 1.3.7]{ramazan}.  

\begin{theorem}
\label{thm:duality}
If $G$ is a second countable locally compact (Hausdorff) Abelian group bundle 
with Haar system then the dual $\widehat{G}$ 
has a dual group bundle, denoted
$\doubledual{G}$. Furthermore, the map $\Phi:G\rightarrow \doubledual{G}$
such that 
$$
\Phi(s)(\omega) = \hat{s}(\omega) := \omega(s)
$$
is a (topological) group bundle isomorphism between $G$ and $\doubledual{G}$. 
\end{theorem}

Before we continue, it will be useful to see that the group bundle notion of duality is a natural extension of the usual Pontryagin dual, 
as illustrated by the following proposition.  
\begin{prop}
\label{prop:3}
Let $G$ be a second countable locally compact Abelian group bundle with 
Haar system.  Then $C^*(G)\cong C_0(\widehat{G})$ via the 
Gelfand transform.  Furthermore, if $f\in C_c(G)$ then the Gelfand transform 
of $f$ restricted to $\widehat{G}_x$ is the Fourier transform of $f|_{G_x}$.  
\end{prop}
\begin{proof}
The first statement follows from the fact that we defined $\widehat{G}$ to be
the spectrum of the Abelian $C^*$-algebra $C^*(G)$.  Next, let $\hat{f}$ be the 
Gelfand transform of $f$.  Then for $\omega\in \widehat{G}$ we see from 
\eqref{eq:1} that 
$$
\hat{f}(\omega) = \omega(f) = \int_{G_{\hat{p}(\omega)}} f(s)\omega(s) d\beta^{\hat{p}(\omega)}(s).
$$
This of course implies that $\hat{f}$ is the usual Fourier transform 
on $\widehat{G}_x$.  
\end{proof}

We can now begin the process of proving Theorem
\ref{thm:duality}.  The first step is to
show that $\widehat{G}$ has a dual bundle.  We have already verified
that $\widehat{G}$ is a 
second countable locally compact Abelian group bundle.  The only
remaining requirement is that $\widehat{G}$ has a Haar system. 
Recall that given a locally
compact Abelian group $H$ and Haar measure $\lambda$ the Plancharel
theorem guarantees the existence of a dual Haar measure
$\hat{\lambda}$ such that $L^2(H,\lambda)\cong
L^2(\widehat{H},\hat{\lambda})$.  The existence of a dual Haar system
is then taken care of by the following lemma.

\begin{lemma}[{\cite[Proposition 3.6]{ctgIII}}]
If $G$ is an Abelian second countable locally compact group bundle
with Haar system $\{\beta^x\}$, then the collection of dual Haar measures 
$\{\hat{\beta}^x\}$ is a Haar system for $\widehat{G}$.
\end{lemma}

Now that $\doubledual{G}$ is well defined, we must show that
$\Phi$ is a group bundle isomorphism.  In some sense, the following proposition gets us most of the way there.  
\begin{prop}
\label{prop:2}
The map $\Phi:G\rightarrow \doubledual{G} : s\mapsto \hat{s}$
is a continuous bijective groupoid homomorphism.  
\end{prop}
\begin{proof}
It follows from Lemma \ref{thm:dual} that $\doubledual{G}_x$ is the double 
dual of $G_x$. 
Furthermore, classical Pontryagin duality says that $s\rightarrow \hat{s}$ is 
an isomorphism from $G_x$ onto $\doubledual{G}_x$ \cite[Theorem 1.7.2]{rudinfourier}.  
Since 
$\Phi$ is formed by gluing all of these fibre isomorphisms together it is 
clear that $\Phi$ is a bijective groupoid homomorphism.  Next, we need to 
see that it is continuous.  Suppose $s_i\rightarrow s_0$ in $G$.  We 
know from Lemma \ref{prop:1} that it will suffice to show that 
\begin{enumerate}
\item $\hat{\hat{p}}(\Phi(s_i))\rightarrow \hat{\hat{p}}(\Phi(s_0))$, and 
\item given $\omega_i\in \widehat{G}_{\hat{\hat{p}}(\Phi(s_i))}$ such that 
$\omega_i\rightarrow \omega_0$ in $\widehat{G}$ then 
$\Phi(s_i)(\omega_i)\rightarrow \Phi(s_0)(\omega_0)$.  
\end{enumerate}

First, let $x_i=p(s_i)=\hat{\hat{p}}(\Phi(s_i))$.  Since $p$ is continuous,
it is clear that $x_i\rightarrow x_0$ and that the first condition is 
satisfied.  Now suppose $\omega_i \in \widehat{G}_{x_i}$ for all $i\geq 0$ 
such that $\omega_i\rightarrow \omega_0$.  All we have to do is 
cite Lemma \ref{prop:1} again to see that 
$$
\Phi(s_i)(\omega_i) = \omega_i(s_i) \rightarrow \omega_0(s_0) = \Phi(s_0)(\omega_0).
$$
\end{proof}

If we were working with groups, we would be done since continuous bijections 
between second countable locally compact groups are automatically 
homeomorphisms \cite[Theorem D.3]{tfb2},\cite[Corollary 2, p. 72]{invitation}.  However,
there currently no automatic continuity results for the inverse of a 
continuous bijective group bundle homomorphism. Regardless, we can
still show that in this case $\Phi$ is a homeomorphism.

\begin{proof}[Proof of Theorem \ref{thm:duality}]
Given Proposition \ref{prop:2}, all we need to do to prove that $\Phi$
is a homeomorphism is show that if $\hat{s_i}\rightarrow \hat{s_0}$ in 
$\doubledual{G}$ then $s_i\rightarrow s_0$ in $G$.  First, we 
let $x_i = p(s_i)$ for all $i$.  Recall that $\doubledual{G}$ has the Gelfand 
topology as the spectrum of $C^*(\widehat{G},\hat{\beta})$.  Therefore, 
for all $\phi\in C_c(\widehat{G})$ we have $\hat{s_i}(\phi)\rightarrow
\hat{s_0}(\phi)$.  When we remember that characters in $\doubledual{G}$
act on functions in $C_c(\widehat{G})$ via equation (\ref{eq:1}) we see 
that this says, for all $\phi\in C_c(\widehat{G})$,
\begin{equation}
\label{eq:2}
\int_{\widehat{G}} \phi(\omega) \omega(s_i) d\hat{\beta}^{x_i}(\omega)
\rightarrow \int_{\widehat{G}} \phi(\omega) \omega(s_0) 
d\hat{\beta}^{x_0}(\omega).
\end{equation}

Now suppose we have a relatively compact open neighborhood $V$ of $x_0$ in 
$G$.  Then using the continuity of multiplication, there exists 
a relatively compact open neighborhood $U$ of $x_0$ in $G$ such that 
$U^2 \subseteq V$.  Choose $h\in C_c(G)$ such that $h(x_0)=1$ and 
$\supp(h)\subseteq U$.  Let $f = h^* * h$.  
Then $f\in C_c(G)$ and a simple calculation shows that $\supp(f)\subseteq V$.  
From now on, let $f^x$ denote the 
restriction of $f$ to $G_x$.  It is clear from the definition of $f$ and 
\cite[Section 1.4.2]{rudinfourier} that it is a positive definite function
on each fibre and 
therefore satisfies the conditions of Bochner's theorem and the inversion
theorem on each fibre.  In particular, it can be shown using
\cite[Section 1.4.3]{rudinfourier} that for each $x$ there exists a finite positive measure $\mu^x$ on $\widehat{G}_x$ 
(extended to $\widehat{G}$ by giving everything else measure zero)
such that 
$$
f(s) = \int_{\widehat{G}} \overline{\omega(s)} \mu^{p(s)}(\omega).
$$
Furthermore, it is easy to prove using \cite[Section 1.4.1]{rudinfourier} that
$\mu^x(\widehat{G}) = \mu^x(\widehat{G}_x) = \|f^x\|_\infty \leq \|f\|_\infty$
for all $x\in X$ so that $\{\mu^x\}$ is a bounded collection of finite 
measures.  
Additionally, it is shown in the proof of \cite[Section 1.5.1]{rudinfourier}
that, as measures on $\widehat{G}_x$, 
$$
\widehat{f^x}d\hat{\beta}^x = d\mu^x.
$$
Proposition \ref{prop:3} states that given $f\in C_c(G)$ the Gelfand transform of $f$ restricts to the 
usual Fourier transform 
fibrewise.  Therefore, since everything outside $\widehat{G}_x$ has measure 
zero, we may as well write 
\begin{equation}
\label{eq:7}
\hat{f}d\hat{\beta}^x = d\mu^x.
\end{equation}

Now, if $\phi\in C_c(\widehat{G})$ then $\phi \hat{f}$ is compactly 
supported.  It follows from (\ref{eq:2}) that
\begin{equation}
\label{eq:8}
\int_{\widehat{G}} \phi(\omega) \hat{f}(\omega) \omega(s_i) d\hat{\beta}^{x_i}
(\omega) \rightarrow \int_{\widehat{G}} \phi(\omega) \hat{f}(\omega) 
\omega(s_0) d\hat{\beta}^{x_0}(\omega).
\end{equation}
Using (\ref{eq:7}), we can rewrite (\ref{eq:8}) as 
\begin{equation}
\label{eq:5}
\int_{\widehat{G}} \phi(\omega)\omega(s_i) d\mu^{x_i}(\omega) \rightarrow 
\int_{\widehat{G}} \phi(\omega)\omega(s_0) d\mu^{x_0}(\omega).
\end{equation}
We can extend (\ref{eq:5}) to functions $\phi\in C_0(\widehat{G})$ by noting 
that $C_c(\widehat{G})$ is uniformly dense in $C_0(\widehat{G})$ and doing a 
straightforward approximation argument using the fact that the 
$\{\mu^{x_i}\}$ are uniformly bounded.  

Let $g\in C_c(G)$. Observe that
\begin{eqnarray*}
\overline{\widehat{g^{x_i}}(\omega)} \omega(s_i) &=&
\int_{G_{x_i}} \overline{g^{x_i}(s)\omega(s)} \omega(s_i) d\beta^{x_i}(s) \\
~&=& \int_{G_{x_i}} \overline{g^{x_i}(s)} \omega(s^{-1}s_i) d\beta^{x_i}(s) \\
~&=& \int_{G_{x_i}} \overline{g^{x_i}(s_i s)} \omega(s^{-1}) d\beta^{x_i}(s) \\
~&=& \overline{(\lt_{s_i^{-1}} g^{x_i})\sidehat}(\omega).
\end{eqnarray*}
Therefore, for all $i$, we have 
\begin{eqnarray*}
\int_{\widehat{G}} \overline{\hat{g}(\omega)} \omega(s_i) d\mu^{x_i}(\omega)
&=&\int_{\widehat{G}} \overline{\hat{g}(\omega)}
\hat{f}(\omega)\omega(s_i)d\hat{\beta}^{x_i}(\omega)\\
~&=& \int_{\widehat{G}_{x_i}} \overline{\widehat{g^{x_i}}(\omega)}
\widehat{f^{x_i}}(\omega)\omega(s_i)d\hat{\beta}^{x_i}(\omega)\\
~&=&\int_{\widehat{G}_{x_i}} \overline{(\lt_{s_i^{-1}} g^{x_i})\sidehat} \widehat{f^{x_i}}d\hat{\beta}^{x_i} \\
~&=& \int_{G_{x_i}} \overline{\lt_{s_i^{-1}} g^{x_i}} f^{x_i} d\beta^{x_i}, 
\end{eqnarray*}
where the last equality follows from the Plancharel theorem \cite[Theorem 1.6.1]{rudinfourier}.  
Since $\bar{\hat{g}}\in C_0(\widehat{G}),$ it follows from (\ref{eq:5}) that
\begin{equation}
\label{eq:3}
\int_{G_{x_i}} \overline{\lt_{s_i^{-1}} g^{x_i}} f^{x_i} d\beta^{x_i} \rightarrow 
\int_{G_{x_0}} \overline{\lt_{s_0^{-1}} g^{x_0}} f^{x_0} d\beta^{x_0}.
\end{equation}

We are now ready to attack the convergence of the $s_i$.  Choose an open neighborhood
$O$ of $s_0$.  Using the continuity of multiplication, we can 
find relatively compact open neighborhoods $V$ and $W$ in $G$ such 
that $x_0\in V$, $s_0\in W$ and 
$VW\subseteq O$.  Furthermore, by intersecting $V$ and $V\inv$ we can 
assume that $V\inv=V$.  Construct $f$ for $V$ as 
in the beginning of the proof. Now choose $g\in C(G)$ so that $0\leq g\leq 1$, 
$g(s_0)=1$, and $g$ is zero off $W$.  
Then $g\in C_c(G)$ and $\overline{g} = g$ so that by equation (\ref{eq:3}) 
we have 
\begin{equation}
\label{eq:4}
\int_{G_{x_i}} g(s_i t) f(t) d\beta^{x_i}(t) \rightarrow
\int_{G_{x_0}} g(s_0 t) f(t) d\beta^{x_0}(t).
\end{equation}
It turns out that $\int g(s_it)f(t)d\beta^{x_i}(t)=0$ unless 
$s_i \in WV^{-1} = WV \subseteq O$.  Furthermore,
both $g(s_0 x_0)$ and $f(x_0)$ are nonzero by construction, and since both 
functions are continuous, this implies
$$
\int_{G_{x_0}} g(s_0 t)f(t) d\beta^{x_0}(t)\ne 0.
$$
It follows from (\ref{eq:4}) that eventually $s_i\in O$.  
This of course implies that $s_i 
\rightarrow s_0$ and we are done.  
\end{proof}

{\em Acknowledgments}:  The author would like to thank the referee for
their helpful comments and for bringing \cite{ramazan} to his
attention.  He also thanks Dana Williams for generously providing 
advice and guidance.

\providecommand{\bysame}{\leavevmode\hbox to3em{\hrulefill}\thinspace}
\providecommand{\MR}{\relax\ifhmode\unskip\space\fi MR }
\providecommand{\MRhref}[2]{%
  \href{http://www.ams.org/mathscinet-getitem?mr=#1}{#2}
}
\providecommand{\href}[2]{#2}

\end{document}